\documentclass[11pt]{article}

\usepackage[centertags]{amsmath}

\usepackage{newlfont,color}

\usepackage{latexsym,enumerate,epsfig,listings}
\usepackage{amsmath,amsthm,amsopn,amstext,amscd,amsfonts,amssymb}
\usepackage{graphicx}

\def \beq {\begin{eqnarray}}
\def \eeq {\end{eqnarray}}
\def \beqn {\begin{eqnarray*}}
\def \eeqn {\end{eqnarray*}}

\newcommand{\beeq}{\begin{equation}}
\newcommand{\eneq}{\end{equation}}
\newcommand{\bearno}{\begin{eqnarray*}}
\newcommand{\enarno}{\end{eqnarray*}}
\newcommand{\befi}{\begin{figure}}
\newcommand{\enfi}{\end{figure}}

\newcommand{\be}{ \begin{equation}}
\newcommand{\ee}{\end{equation}}

\def\ep{\hfill $\Box$}

\def\bp{\noindent{\bf Proof.}\ }

\def\P{\mathbb P}
\def\R{\mathbb R}

\def\C{{\cal C}}

\def\U{{\cal U}}

\def\A{{\cal A}}

\def\S{{\cal S}}
\def\P{{\cal P}}

\newtheorem{lemma}{Lemma}[section]

\newtheorem{remark}{Remark}[section]
\newtheorem{definition}{Definition}
\newtheorem{example}{Example}
\newtheorem{theorem}{Theorem}
\newtheorem{proposition}{Proposition}

\title{Large deviations of the stationary measure of networks under proportional fair allocations}

\author{M. Jonckheere, S. L\'opez}

\begin{document}

\maketitle

\begin{abstract}

We address a conjecture introduced by Massouli\'e (2007), concerning the large deviations of the stationary measure of bandwidth-sharing
networks functioning under the Proportional fair allocation. For Markovian networks, we prove that Proportional fair and an associated reversible allocation are geometrically
ergodic and have the same large deviations characteristics using Lyapunov functions and martingale arguments. For monotone networks, we  give a more direct proof of
the same result relying on stochastic comparisons that hold for general service requirement distribution. These results comfort the intuition that Proportional fairness is `close' to
allocations of service being insensitive to the service time requirement.
\end{abstract}

\section{Introduction}

Bandwidth-sharing networks describe the evolution of the number of flows (or calls) in a communication network where different classes of traffic compete for the bandwidth.\\

They have become a standard modeling tool over the past decades for modeling communication networks~\cite{bp,massrob} and have been used in particular
 to represent the flow level dynamics of a wide range of wireline and wireless
networks \cite{BFreview}, 
generalizing henceforth more traditional voice traffic models~\cite{kelly1979}. \\

In queuing theory, these models boil down to a particular class of processor sharing networks with state-dependent service rates. Assuming that class-$i$ flows arrive subject to a Poisson
process of intensity $\lambda_i$ and require exponentially
distributed service times of mean $\mu_i$ (the arrival processes of all classes being mutually
independent), and in the absence of internal routing, the stochastic process $X=(X_1, \ldots, X_N)$
describing the number of flows (or calls) in progress in the network is a multi-dimensional birth and death process with transition rates:
\begin{eqnarray*}
q(x,x-e_i)&=&\mu_i \phi_i(x),\\
q(x,x+e_i)&=&\lambda_i,
\end{eqnarray*}
where $x=(x_1,\ldots,x_N)$ is the number of flows in each class of traffic.
The service rates of the $N$ traffic classes $\phi=(\phi_i(\cdot))_{i=1}^N$ encodes the particularities of the network
resulting from the specific topology, technology, radio conditions, interference and multi-diversity effects and the packet protocols  and congestion control mechanisms in use. In a wireline network, the vector $\phi(x)$ is usually assumed to belong for all $x$ to a polyhedron describing the capacities constraints of each link that are used by the different routes while for wireless networks,  $\phi(x)$ generally belongs to a more complicated closed convex set containing $0$
 corresponding to the achievable rates.
 In many situations like TDMA or CDMA data networks, the capacity set can be assumed to be convex but this is not necessarily the case
for decentralized schemes like for the 802.11 WLAN. However, it has been shown that the capacity sets in that case are log-convex (see \cite{duffy2010}).

\

 Some specific bandwidth allocations have received a lot of attention in recent years.
These include the max-min fair allocation and the proportional
fair allocation (PF) that maximize a $\log$-utility function. More generally,
Mo and Walrand introduce the following family of utility functions:
$$U^\alpha(x,\eta)= \sum_{i=1}^N x_i{\eta_i^{1-\alpha} \over 1-\alpha}, \qquad \alpha \in (1,\infty),$$
$$U^1(x,\eta)=  \sum_{i=1}^N x_i \log(\eta_i) , \qquad \alpha =1,$$
including as special cases the proportional fair allocation and the max-min fair allocation for $\alpha$ tending $\infty$.
On the other hand, the balanced fairness allocation (BF) was defined in \cite{perf}  as the  allocation ensuring the reversibility of the Markov process $X$
and maximizing the probability of the network to be empty (among the `reversible' allocations). \\

More generally, the selection of a specific bandwidth allocation
is motivated by several properties of the resulting process $X$.
Among those properties, the following ones are of particular interest:
\begin{enumerate}
\item
{\it maximal stability}: one expects the bandwidth sharing mechanism to stabilize the system whenever it can be stabilized,

\item
{\it decentralized protocol}:
the bandwidth allocation can be implemented in the network using decentralized schemes,

\item
{ \it robustness}: when changing the traffic conditions\footnote{When the size distributions of the flows are not exponentially distributed, the process $X$ is not Markov by itself anymore and the dynamics have to be defined using the residual service time.} and in particular the service time distribution (but keeping the mean flow size fixed), one could expect the
stationary measure of $X$ to remain the same, in which case the system and the bandwidth allocation are said to be insensitive,

\item
{ \it relaxed robustness }: the large deviations characteristics do not depend on the service time distribution, except for its mean and coincide with the large deviations of the most efficient insensitive allocation.

\end{enumerate}

The PF allocation satisfies properties (1) and (2) but fails to satisfy property (3) on general topologies. On the other hand, the BF allocation satisfies (1) and (3) but it is not known whether it satisfies (2) \cite{BFreview}.
On the one hand, it was shown in \cite{massoulie2007} that an appropriate modification of the proportional fair allocation, called modified proportional fair allocation,  (mPF)
 coinciding asymptotically (point-wise) with PF, has the same large deviations characteristics as BF.
On the other hand, it has been recently proven that an insensitive allocation being maximal stable,
is asymptotically equivalent (point-wise) to PF \cite{walton2011}.
However, it remained an open problem to prove that the large deviation characteristics of the stationary measure of the PF allocation itself coincide with those of mPF,
and BF as it was conjectured in \cite{massoulie2007}. \\

Our contribution is the following.
We first recall the principles of the Freidlin and Wentzell theory for birth and death processes on $\mathbb Z^N$ with rates being Lipschitz and with bounded logarithms.
As underlined for instance in \cite{schwartzweiss}, it is very demanding to extend these results to processes on state spaces
with boundaries since the technical conditions of the classical theory are never fullfilled. (Remark that  boundness of the logarithm is never verified for birth and death processes living in the orthant,
 while the Lipschitz assumption is not verified in our case for several network topologies). \\

 A main contribution of the present article consists in overcoming these difficulties for the specific processes we are studying.
We first show that for Poisson arrivals and exponentially distributed flow sizes, the stationary distribution $\pi^{PF}$ of the number of flows associated with PF and the stationary 
distribution $\pi^{mPF}$ associated with the mPF allocation
have the same large deviations characteristics. More precisely:
\begin{theorem}
For all $x \in \mathbb Z _{+}^N$:
$$ {1 \over n }\log \Big( {\pi^{PF}(n x) \over \pi^{mPF}( nx  )} \Big) \leq O(n^{-\frac{1}{2}+ \epsilon}), \qquad \forall \epsilon >0.$$
In the particular case that $x \in \mathbb N ^N$ the bound can be improved to  $O(n^{-1})$.
\end{theorem}
This is achieved by first proving the {  geometric ergodicity} of both the mPF allocation and the PF allocation. For that purpose, we exhibit appropriate Lyapunov functions,
 relying on some structural results of PF described in \cite{massoulie2007}.
This then allows us to use simple martingale arguments.

\

Finally, for { monotone} networks and generally distributed flow size, we give a more direct proof establishing that the large deviations characteristics are actually insensitive
 to the service time distribution. This shows that the proportional fair allocations indeed satisfy properties (1) and (2) and (4) at least on monotone topologies.

\

The rest of the paper is organized as follows.
In Section \ref{sec:allocations}, we introduce the  allocations functions used in the sequel as well as several of their properties. In Section \ref{sec:FW}, we show that the rate function
of the proportional fair allocation coincides for Markovian dynamics to the rate function of Balanced fairness. In Section \ref{sec:MN}, we prove property $(4)$ for monotone networks using stochastic comparisons.

\section{Properties of bandwidth sharing allocations}\label{sec:allocations}

\paragraph{Notations}

We define here a few notations that we need in the sequel.
$\mathbb R_+$ denotes the set of non-negative real numbers.
The bandwidth allocation vector is denoted $\phi(x)\equiv (\phi_i(x))_{i=1 \ldots N}$.
For any vector $v$ in $\mathbb{R}^N$, and function $f: \mathbb R \to \mathbb R$, we denote $f(v)$ the vector $(f(v_i))_{i=1 \ldots N}$.
Similarly, we  also use the notation $a^v$ for $\prod_{i=1}^N a_i^{v_i}$.
The usual scalar product in $\mathbb R^N$ between $u$ and $v$ is denoted $\langle u,v \rangle$. We use $|| \cdot ||_p$ to denote the $l_p$-norm, but we reserve
$|\cdot |$ to denote the $l_1$-norm: $|v|=\sum_{i=1}^N |v_i|$.
For $x,y \in \mathbb R^N$, we also use the notation $x \le y$ to denote the partial order $x_i \le y_i$ for all $i=1 \ldots N$.
Finally, given a set $\C$ and function $f: \mathbb R \to \mathbb R$, we shall denote $f(\C)$ the set $\{ f(\eta) : \eta \in\C \}$.

\subsection{Gradient allocations}

In the characterization of the large deviations of the stationary regime, we shall rely on
 two properties of the allocation function $\phi$ playing a crucial role and being intrinsically related:
being a gradient allocation or a discrete gradient allocation.

\begin{definition}
A bandwidth allocation is called gradient if there exists a function $P : \R^N \to \R $ (that we call a continuous potential) such that
$$\log(\phi(x))= -\nabla P(x), \quad \forall  x \in \mathbb R_+^N \setminus\{0\}.$$
\end{definition}

Recall that the proportional fair allocation is defined by the optimization problem associated to a capacity set $\C$ as follows:
$$\phi^{PF}(x)=\arg \max_{\eta \in \C}U^1(\eta,x)=\arg \max_{\eta \in \C} \langle x,\log(\eta) \rangle.$$

Following Massouli\'e \cite{massoulie2007}, observe that the proportional fair bandwidth allocation is gradient.
Indeed, let $\delta^*_{\A}$ the support function of a bounded convex set $\A$ i.e.:
$$\delta^*_{\A}(x)=\max_{\eta \in \A}{ \langle x,\eta \rangle}.$$

\begin{proposition}
Assume that the set $\C$ is log-convex, i.e., the set $\log(\C)$ is convex,
then
$$\log(\phi^{PF}(x))=  \nabla \delta^*_{\log(\C)}(x), \quad \forall  x \in \mathbb R_+^N \setminus \{0\}.$$
\end{proposition}
\bp

The function $\delta^*$ is sub-differentiable because it is convex and finite (see \cite{rockafellar1970})  for all  $x \in \mathbb R_*^N$.
The unicity of the sub-gradient comes from the strict concavity of the $\log$ function and implies the differentiability.

\ep

In the sequel, we always assume that $\C$ is convex and contains the set $\{ \eta: \sum \eta_i \le c\}$ for some $c>0$.
Of course, this is not a restriction for applications. We denote by $P^{PF}\equiv \delta^*_{\log(\C)}$.

\subsection{Reversible allocations}

\begin{definition}
A bandwidth allocation is called  reversible or discrete gradient,  if there exists a function $\tilde P : \R^N \to \R $ (that we call a {\it discrete potential})
such that for all $x \in \mathbb N^N$, $i=1,\dots, N$:
$$\log(\phi_i(x))= - \mathbf{D} \tilde P(x) \equiv \tilde P(x)- \tilde P(x-e_i).$$
\end{definition}

In the latter case, the stationary measure of the process is easily described:

\begin{proposition}\label{prop:rev}
The process $\tilde X$ associated with a reversible allocation with discrete potential $\tilde P$ is reversible (in the usual sense) and its stationary measure (when it exists) is:
$$\pi(x)= C \lambda^x \exp(\tilde P(x)).$$
Assume further that ${1\over n} \tilde P(nx) \to \gamma(x)$ as $n \to \infty$.
Then:
$$ \lim_{n \to \infty}{1\over n}\log\pi(n x)=-\big(\gamma(x)-\sum_{i=1}^N x_i \log(\lambda_i) \big).$$
\end{proposition}

A particular role in what follows shall be played by 
two reversible allocations `close' to PF: 
 the modified  PF allocation (mPF) and the BF allocation which we now define mathematically.

\begin{definition}[Modified proportional fair and Balanced fair allocations]
$$\log(\phi^{mPF}(x))= - \mathbf{D} \tilde P^{PF}(x),$$
with $P^{PF}$ the continuous potential associated with PF, and
$$\log(\phi^{BF}(x))= - \mathbf{D} \tilde P(x),$$
where the potential $\tilde P$ is recursively defined by:
$$\tilde P(x) =0,$$
$$ \tilde P(x) = \max \{a>0 :   \tilde P(x-e_i) - a \in \log(\C) \} .$$
\end{definition}

We have the following large deviations results for the BF allocation:

\begin{proposition}{ (Massouli\' e, 2007)}\label{prop:BF_LD}
$$ \lim_{n \to \infty}{1\over n}\log\pi^{BF}(n x)= P^{mPF}(x).$$
\end{proposition}

The main challenge of the subsequent analysis is to provide similar results for
the Proportional fair allocation itself and hence verify that both allocations share the same large deviations characteristics.
The main difficulty of this program is that the stationary measure of a network under the proportional fair allocation does not have a closed form in general (except in the very particular cases of symmetric hypergrids).

\subsection{Stability}

We state two stability results, one for the original stochastic system, and one for the deterministic analogue %(that can be proved to be the fluid limit, see the next Section).
Both results are well known for Markovian bandwidth sharing networks in \cite{PEreview}.
Massouli\'e extended the stability results to general service time distributions in \cite{massoulie2007}.

\begin{theorem} [\cite{massoulie2007}]
The stability set of the process $X$ associated with the PF allocation
 is contained in, and contains the interior of the set $\S=\bar \C.$
\end{theorem}

\begin{remark}
For Markov processes, a direct proof follows from taking $x \mapsto \sum_{i=1}^N {x_i^2 \over \lambda_i}$ as a Lyapunov function.
\end{remark}

Define now the deterministic dynamical system  $(DS)_{x_0}$ by
\begin{eqnarray}
&\dot{x_i}&=\big(\lambda_i-\phi_i(x)\big)1_{x_i>0},\\
&x(0)&=x_0.
\end{eqnarray}

It was proved \cite{kellywilliams} (see also \cite{massoulie2007}) that the fluid limits
of the process $X$ are solutions of this dynamical system at all regular points.

\subsection{Examples}

We illustrate the above description on a few toy examples.

\begin{example}[Single link]
Consider first a single link shared by $N$ classes.
This corresponds to choosing $\C=\{ \sum {\eta_i \over \mu_i} \le 1 \}$, where the $\{ \mu_i \}$
correspond to the mean flow sizes.
In this case, the  proportional fair allocation coincides with the balanced fairness allocation.
Its service rates  are given by:
$$\phi_i(x)=\mu_i{x_i \over |x|}.$$
Note however that the continuous and the discrete potential functions associated with the allocation $\phi$ do not coincide:
$$\delta^*_{\log(\C)}(x)=\sum_{i=1}^{N} {x_i \over |x|} \log \left( {x_i \over |x|}\right),$$
while the discrete potential is given by:
$$ P(x)= \log \left({|x| \choose x_1,\ldots, x_N} \right).$$
One can however prove using the Stirling formula that
$${1 \over n} (P(n x)-\delta^*_{\log(\C)}(n x)) \to 0, \textrm{ as } n \to  \infty.$$
An extension of these formulae can be obtained for hypergrids topologies \cite{bp}.
\end{example}

\begin{example}[Tree network]
A 2-level tree network with $N$ routes is defined by the following polyhedron for the capacity set:
$$\P=\{\eta  \in \mathbb R_+^N:   \eta_i \le c_i, ~ \sum_i \eta_i \le 1  \    \},$$
for some constants $\{c_i\}$.
An expression of the stationary measure for the balanced fairness allocation
can be found in \cite{bp}. In general (i.e., when the tree is not degenerate), the proportional fair allocation does not coincide with the balanced fair allocation and does not
raise a reversible process. Hence its stationary distribution can be obtained only numerically.
For the case $N=2$,  and $\max(c_1,c_2)\le 1$, the proportional fair allocation is:
$$\phi_1(x)=1-\phi_2(x),$$
$$\phi_2(x)=\min(c_2,{x_2 \over x_1+x_2}).$$

An important property of the proportional fair allocation in that case is its monotonicity (see Section 4 for a definition) \cite{stobounds}:
$${\phi_i^{PF}(x) \over x_i} \ge {\phi_i^{PF}(y) \over y_i},  ~ x_j  \le  y_j, ~ \forall j. $$

\end{example}

\begin{example}[Wireless network]
In \cite{duffy2010}, the rate region $\C$ of a $2$-station network functioning under the 802.11e protocol is studied.
The rate region is the set of achievable throughput vectors
at the fine time scale (packet level). Their findings show that the rate region, (which exact expression depends in a complicated manner of the probabilities of transmitting) is generally not convex but is however $\log$-convex.
\end{example}

\pagebreak

\section{Large deviations for Markovian dynamics}\label{sec:FW}
\subsection{The Freidlin-Wentzell theory for birth and death processes with smooth transitions}

In this Section, we remind the concepts introduced by Freidlin and Wentzell \cite{freidlin1984}
allowing to get a grasp on the sample-paths large deviations of stochastic
differential equations with small noise. These results were often applied to diffusions processes
but can equally be applied (as indicated in \cite{freidlin1984}) to { some} jump processes.
We show here how these classical results can be used to prove a large deviation principle for birth and death processes  with ``gradient'' rates on $\mathbb Z^N$. \\

 Let $Y^n$  be a multi-dimensional birth and death process with transition rates:
\begin{eqnarray*}
q\left({x \over n},{x \over n}-{e_i \over n}\right)&=&n  \phi_i\left({x \over n}\right),\\
q\left({x \over n},{x \over n}+{e_i \over n}\right)&=& n \lambda_i.
\end{eqnarray*}
We suppose additionally that the death rates are  $0$-homogeneous (i.e. $\phi(a z )= \phi(z)$, $\forall a >0$). \\

Define the logaritmic moment-generating of the increment of the process for $z \in \mathbb R^N/n, y \in \mathbb R^N$ by
\beeq
H_n({z},y) ={d \over dt} E^{z}[\exp( \langle y,Y^n(t) -{z} \rangle)].
\eneq
Using the structure of the generator and the $0$-homogeneity of the rates:
\beeq
H_n({z},y)= n \left( \sum_{i} \lambda_i  (e^{y_i/n}-1) + \phi_i({z})(e^{-y_i/n}-1)   \right).
\eneq
We thus obtain that $H_n({z},y/n) = n H(z,y/n)$ with
\beeq
H(z,y)= \langle e^y-1,\lambda \rangle + \langle e^{-y}-1,\phi(z) \rangle.
\eneq

Define now the Fenchel-Legendre $L$ transform of $H$:
\beeq
L(x,y)=\sup_{\theta \in \mathbb R^N} \langle y,\theta \rangle -H(x,\theta),
\eneq
and define the action functional $S: \C([0,T]) \to \mathbb R$ as follows:
\beeq
S_T(r)  = \int_{0}^T L(r_s,\dot{r}_s)ds.
\eneq
Now define the quasi-potential $V$ by:
 \beeq
V(x)=\min_{r, T, r(0)=0, r(T)=x} S_T(r,\dot{r}).
\eneq

Assume that the function $L$ is such that:
$$\sup_{|x-\tilde x| } {|L(x,y)-L(\tilde x,y)| \over 1+ L(x,y)} \to 0.$$

The original results of Freidlin and Wentzell are stated for Lipschitz-continuous transitions (extended to $\mathbb R^N$):

\begin{theorem}[Freidlin and Wentzell]
Assume that for each $i$, the function $\log \phi_i(\cdot)$ is bounded and Lipschitz continuous.
Let $\pi^n$ the stationary measure of the birth and death processes $X^n$.
Then, a large deviations principle holds for the family of probabilities $\pi^n(\cdot)$ with rate function $V(\cdot)$.
\end{theorem}

\subsubsection{Further characterization of the potential}

In the case of general multi-dimensional birth and death processes, it is difficult to solve the variational problem from which the potential $V$ is defined.
However, assuming that the rates are gradient, as defined in the previous Section, greatly simplifies the expression of the potential.

\begin{proposition}\label{prop:pot}
Assume that the rates (allocation) $\phi$ is gradient with continuous potential $P$.
Then the quasi-potential of the birth and death process is equal to $R(x)=- \langle \log(\lambda),x \rangle +P(x)$.
\end{proposition}

\bp

We first show that:
$$V(x)=\max_{x_s : x_0=0, x_T=x} \int_0^T \langle {\dot x}_s ,\log \Big( {\lambda \over \phi(x_s)} \Big) \rangle ds.$$
The first simplification comes from the time homogeneity of the process.
Using the Euler-Lagrange principle,
combined with the time homogeneity (which implies that $\nabla_x {L}=0$)  the minimizing path satisfies the Beltrami identity:
$$L- \langle y, \nabla_y{L} \rangle=0.$$
Now, observe that if $\theta_y =\arg\max \langle y,\theta \rangle -H(x,\theta)$, we obtain that
$ \nabla_y{L}=\theta_y$ and $L-\langle y \nabla_y{L} \rangle =H(x,\theta_y)=0.$
Solving the last equation leads to:
$$\theta_{\dot x}={1 \over 2} \log({\phi(x) \over \lambda}),$$
which allows us to conclude since
$$L(x,\dot x)= \langle \dot x, \nabla_y{L(x,\dot x)} \rangle = \langle \dot x, \theta_{\dot x} \rangle.$$

 The expression for the quasi-potential then follows from integration.

\ep

We insist on the facts that this theory
allows to prove the large deviation rate of the stationary measure of birth and death processes
only for smooth allocations and that it becomes a very technical issue to weaken this assumption. Also it does not provide a rate of convergence in general.

We address these issues in the next Section: we prove directly the large deviation principle
using martingale arguments by taking advantage of the geometric ergodicity of the process.

\subsection{Large-deviations results for general gradient allocations}

In this Section, we call $X$ the process associated with the PF allocation, denoted $\phi$ and corresponding to the { \bf continuous gradient} of the potential $P^{PF}$.
We further denote by $\tilde X$ the process associated to the allocation $\tilde \phi$
itself corresponding to the {\bf discrete gradient} of $P^{PF}$, i.e.:
\begin{eqnarray}
\log(\phi(x))&= &\nabla P^{PF}(x),\\
\log(\tilde \phi(x) ) &=& -\mathbf{D} P^{PF}(x) = ( P^{PF}(x)-P^{PF}(x-e_i))_{i=1\ldots N}.
\end{eqnarray}

\subsubsection*{Structural properties of the PF allocation}

From the structural representation of PF, we know that $\tilde{\phi}$ is a perturbation of the original rates $\phi$ in the sense that (see Lemma 9 in \cite{massoulie2007} ):
\beeq\label{def.pert}
\big|\log( \phi_i(x))-\log(\tilde \phi_i(x))\big| = \big| \log(\phi_i(x))- (P^{PF}(x)-P^{PF}(x-e_i))\big|\le  {1\over x_i} , ~ \forall i, x.
\eneq

\pagebreak

\subsubsection*{Geometric ergodicity}

We first consider the dynamics corresponding to the Proportional fairness allocation.

\begin{proposition}[Geometric ergodicity of $ X$]
Suppose $\lambda \in int(\C)$, then there exists a constant $K$ such that for $|x|>K$:
$$\Delta G(x) \le - \gamma G(x).$$
Hence:
$$\Big| P^0_t(x) - \pi^{PF} (x) \Big| \leq K_1 e^{ - K_2 t },$$
for some constants $K_1,K_2>0$, and $\pi^{PF}$ the stationary distribution of ${X}$.
\end{proposition}

\bp

The proof has two steps. We first construct a Lyapunov function with bounded drift.
We then use this Lyapunov function to construct a new one verifying a geometric drift inequality.

\paragraph{First step:} let $F(x)=(\sum_{i=1}^N {x_i^2 \over \lambda_i})^{1/2}$.
Observe that $F$ is a norm   in $\mathbb R^N$ (hence is  positive, $1-$ homogeneous, and diverges to infinity when $|x| \to \infty$).
Furthermore, it is $C^2$ for all $x \neq 0$ and:
$$\nabla F(x)=  {x_i \over \lambda_i F(x)}.$$
Hence, there exists $K>0$ such that for $|x| > K$
$${1\over |x|}\sup_{z=x, x \pm e_i }|{\partial^2 F(z) \over \partial^2 x_i}|  \le \epsilon.$$
Using that $F$ is $1$-homogeneous, this leads for $|x|>K$ to:  
\begin{eqnarray*}
\Delta F(x) & \equiv & \sum_{i=1}^N \lambda_i ({F(x+e_i)}-{F(x)} )+\phi_i(x)({F(x-e_i)}-{F(x)}  ), \\
 & \le &    \langle \lambda  - \phi(x) , \nabla F(x) \rangle + \epsilon, \\
 & = &    \langle \lambda  - \phi(x) , {x \over F(x)\lambda} \rangle + \epsilon, \\
 & = &  {|x| \over F(x)}  \langle \lambda  - \phi(x) , {x \over |x| \lambda} \rangle + \epsilon.
\end{eqnarray*}
Using the definition of Proportional fairness (as the allocation maximizing $U^1(\cdot,x)$ and the strict concavity of the $\log$ function), for all $\eta \in \C$,
$$ \sum_{i=1}^N{x_i \over \eta_i} (\eta_i -\phi_i(x)) = \langle \big({\partial U(x,\eta) \over \eta_i}\big) ,\eta-\phi \rangle <0,$$
which implies that that there exists $\gamma$ such that:
$$ \langle \lambda  - \phi(x) , \frac{x}{|x|} {1\over \lambda} \rangle \le - \gamma,$$
which in turn implies (together with the fact that $F$ is norm-like)
 that there exists $\tilde \gamma$ such that:
\begin{eqnarray*}
\Delta F(x) & \le &  - \tilde \gamma, \quad  \forall x, |x|>K .
\end{eqnarray*}
Remark also that
\begin{eqnarray*}
|\Delta F(x)| & \le &  C,  \quad  \forall x .
\end{eqnarray*}

\paragraph{Second step:} we now calculate the drift of $G(x)=\exp(\delta F(x))$:
\begin{eqnarray*}
\Delta G(x) & = & \sum_{i=1}^N \lambda_i (e^{\delta F(x+e_i)}-e^{F(x)}) )+\phi_i(x)(e^{\delta F(x-e_i)}-e^{\delta F(x)})  ), \\
 & = & G(x)  \sum_{i=1}^N \lambda_i (e^{\delta( F(x+e_i)-F(x)}) -1 )+\phi_i(x)(e^{\delta(F(x-e_i)-F(x)} -1)  ), \\
 & \le & G(x) (\delta \Delta F(x) +  C_1 \delta^2 \sum_{i=1}^N \exp( c |\Delta F(x)| )  + O( \delta^3) ), \\
& \le &  G(x) (- C_2 \delta  + C_3 \delta^2),\\
& \le &  -\gamma   G(x).
\end{eqnarray*}

This implies that (see \cite{meyntweedie1993}):
$$|P^0(X(t)=x)-\pi^{PF} (x)| \le K_3 \exp^{ - K_4 t }.$$

\ep

We now deal with the dynamics of the modified proportional fair allocation.
\begin{proposition}[Geometric ergodicity of $\tilde X$]
Suppose $\lambda \in int(\C)$, then there exists a constant $\tilde K$ such that for $|x|> \tilde K$:
$$\Delta G(x) \le - \gamma G(x).$$
Hence,
$$\Big| P^0_t(x) - \pi^{PF}(x) \Big| \leq K_3 e^{ - K_4 t },$$
for some constants $K_3,K_4>0$, $\pi^{mPF}$ being the stationary distribution of $\tilde X$.
\end{proposition}

\bp

We proceed as previously. 
The only difference consists in proving that the perturbations of the rates are small enough
to be negligible in the drift calculations.
Let $F(x)=(\sum_{i=1}^N {x_i^2 \over \lambda_i})^{1/2}$.
For $|x|> \tilde K$, using the bounds on the difference of $\phi$ and $\tilde \phi$, we obtain that:
$$ \big|\phi_i(x)-\tilde \phi_i(x)\big| \le c_0 \big|\log( \phi_i(x))-\log(\tilde \phi_i(x))\big| \le {c_0 \over x_i}.$$
Hence:
\begin{eqnarray*}
\Delta F(x) & = & \sum_{i=1}^N \lambda_i ({F(x+e_i)}-{F(x)} )+\tilde \phi_i(x)({F(x-e_i)}-{F(x)}  ), \\
 & \le &   \langle \lambda  - \tilde \phi(x) , \nabla F(x) \rangle + \epsilon, \\
 & \le &    \langle \lambda  -  \phi(x) , {x \over F(x)\lambda} \rangle + \langle \phi(x)-\tilde \phi(x) , {x \over F(x)\lambda} \rangle + \epsilon, \\
& \le &   - \tilde \gamma + \langle \phi(x)-\tilde \phi(x) , {x \over F(x)\lambda} \rangle ,\\
& \le &   - \tilde \gamma +  {c \over |x|}.
\end{eqnarray*}
which gives (together with the fact that $F$ is norm-like)
 that there exists $\hat \gamma$ such that:
\begin{eqnarray*}
\tilde \Delta F(x) & \le &  - \hat \gamma, \quad \forall x, |x|> \tilde K .
\end{eqnarray*}
Remark also that
\begin{eqnarray*}
|\tilde \Delta F(x)| & \le &  C .
\end{eqnarray*}
We can conclude following exactly as in the previous proof.\ep \\

Before proving the main result of this section, we establish a useful lemma using the geometric ergodicity of the processes.

\begin{lemma}\label{lem:contpi}
 There exists $C>0$ such that for $t_n= C n$, we have
$$ \log \frac{ \mathbb{P}^0 ({X}_{t_n} \geq nx) }{ \pi^{PF}( \{ nx \}^{\uparrow} )} \leq  C_1 \exp \{ - C_2 \, n \} $$
for $C_1,C_2 >0$, where $\pi^{PF} ( \{ k \}^{\uparrow} ) = \sum_{j \geq k} \pi ^{PF}(j)$. The analogue inequality is true for process $\tilde X$.
\end{lemma}

\bp 

We need to bound from below the stationary probabilities of the process.
This is easy since the rates are bounded. Let us define a process with the same arrival rates and the service rates equal to $\bar \phi$, the maximum  of $\phi_i$ for all coordinates.
It is clear that
$$\pi^{PF} (\{nx\}^{\uparrow}) \geq K_1 \Big({\lambda \over \bar \phi}\Big)^n,$$
for some constant $K_1>0$. \\

We have
\begin{eqnarray*}
\log \frac{ \mathbb{P}^0 ({X}_{t_n} \geq nx) }{ \pi^{PF}( \{ nx \}^{\uparrow} )}  & \leq & \Big| 1 - \frac{ \mathbb{P}^0 ({X}_{t_n} \geq nx ) }{ \pi^{PF}( \{nx \}^\uparrow )} \Big| \\
& = & \Big| \frac{ \mathbb{P}^0 ({X}_{t_n} \geq nx ) - \pi^{PF}( \{nx \}^\uparrow ) }{\pi^{PF}( \{nx \}^\uparrow )} \Big| \\
& \leq & \frac{ K_2 e^{ - K_3 t_n } }{ \pi^{PF}( \{nx \}^\uparrow ) }  \\
& \leq & \frac{ K_2 e^{ - K_3 t_n } }{ K_1 \big({\lambda \over \bar \phi}\big)^n } = \frac{K_2}{K_1} \exp \{ -K_3 \, t_n + \log \Big( {\bar \phi \over \lambda } \Big) n  \}
\end{eqnarray*}
where first inequality comes by $ \log (x) \leq |1-x|$ for all $x>0$, and the second one by the ergodicity of ${X}$. We take $C>0$ such that $K_3 \, C > \log
\big( {\bar \phi \over \lambda } \big)$. Note that the same argument functions for process $\tilde X$ taking its respective transition rates and stationary distribution.
\ep 

\subsubsection*{Change of measure and control of the martingale}

To relate the distribution of $X$ and $\tilde X$, we recall the Proposition B.6 513 of Schwartz and Weiss \cite{schwartzweiss}:

\begin{proposition}
Let $Y^i_t$ for $i=1,2$ two multidimensional birth-death process in $\mathbb{Z}_+ ^N$ with step directions $e_j$, bounded step rates $\{ q_j^i(x) \}_{ x \in \mathbb{Z}_+^N}$ and
law $P^i$. Assume that for all $x$ and $j$, $q_j^1(x)=0$ if and only if $q_j^2(x)=0$. Then we can relate the distributions of the processes by $dP^2 = M_t d P^1$ where
$$ M_t =  \exp \Big( \int_0^t \sum_{j=1}^N  \, \log \Big( { q^2_j ( X_{s^-} ) \over q^1_j ( X_{s^-} ) } \Big) dN^j_s \,
  -   \int_0^t \sum_{j=1}^N  ( q_j^2 (X_{s}) - q_j^1 (X_{s}) ) ds \Big). $$
Here $\{ N^j \}$ denote a family of counting processes describing the jumps of the processes in the $j$-th direction
and $M_t$ is a c\`adl\`ag martingale.
\end{proposition}

In our case, $\phi_i(x)$, $\tilde \phi_i(x)$ are both positive whenever $x_i>0$, and $0$ otherwise. All the rates of $X$ and $\tilde X$ are bounded so we meet the conditions of
the previous Proposition and we can write that $d\tilde P(\omega)= M_{t} d P(\omega),$ with
\begin{eqnarray*}
M_{t} & = & \exp \Big( \int_0^t \sum_{j=1}^N 1_{ \{ X{s^-} >0 \}} \, \log \Big( {\tilde \phi_j(X_{s^-}) \over \phi_j(X_{s^-})} \Big) \, dN^j_s \\
 &  & \ \ \ -  \int_0^t  1_{ \{ X_{s} >0 \} } ( \tilde \phi_j(X_{s}) - \phi_j(X_{s}) ) ds \Big),\\
&=& \prod_{j=1}^N M^j_t,
\end{eqnarray*}
and each of the $M^j$ is itself a martingale.

\begin{remark}
An important observation for the following is that the counting processes $N^j$ can be seen as thinning (according to $X_s$) of some Poisson processes
$ \hat{N}^j$ which are all independent.   
\end{remark}

We denote by $M_{s,t}$ the last expression with the integral running from $s$ to $t$, so $M_t=M_{0,t}$. \\

\noindent The change of measure formula is then:
\beeq
E[ 1_{ \{ \tilde X_t =nx \} }]= E[ M_t 1_{ \{  X_t=nx \} } ].
\eneq

\noindent We can now prove our main result. 
\\

\textbf{Proof of Theorem 1:} \\ 
We consider first the case where $x_i>0$ for all $i=1\ldots,N$.

Define the sequence of stopping times:
\begin{eqnarray*}
\tau_1&=&0,\\
\tau_i &=& \inf\{ t > \tau_ {i-1} , \tilde{X}_t \ge  n x\}, \text{ for } i \text{ even, }  i \ge 2, \\
 \tau_i &=& \inf\{ t > \tau_ {i-1}, \tilde{X}_t  <  n x\}, \text{ for } i \text{ odd, }  i \ge 3
 \end{eqnarray*}

Observe that if $\tilde{X}_t \ge nx$, then necessarily, there exists $k$ even (a.s. finite) such that $\tau_k \le t \le \tau_{k+1}$.
Using the Markov property and the martingale property, we then have
\begin{eqnarray*}
 E^0 ( M_{0,t}  1_{ \{\tilde X_{t} \ge nx \}} )
&=& E^0 \sum_{k \, even}  1_{ \tau_k \le t < \tau_{k+1} } M_{t}  1_{ \{\tilde X_{t} \ge nx \}}\\
 &=& \sum_{k \, even} E^0  1_{\tau_k \le t <\tau_{k+1} }  M_{0,\tau_k} M_{\tau_k,t} 1_{ \{\tilde X_{t} \ge nx \}},\\
  &=& \sum_{k \, even} E^0 ( E ( 1_{\tau_k \le t <\tau_{k+1} }  M_{0, \tau_k} M_{\tau_k,t} 1_{ \{\tilde X_{t} \ge nx \}} | \tilde{\mathcal{F}} _{\tau _k} ) ) \\
 &=& \sum_{k \, even}  E^0 ( M_{0,\tau_k} E ( 1_{\tau_k \le t <\tau_{k+1} }  M_{\tau_k,t}1_{ \{\tilde X_{t} \ge nx \}} | \tilde X_{\tau _k} ) ) ,\\
\end{eqnarray*}
with $\{ \tilde{\mathcal{F}} _{k} \}$ being the natural filtration of process $\tilde X$. \\

We define
$$ g_k= \sup_{y} E ( 1_{\tau_k \le t <\tau_{k+1} }  M_{\tau_k,t}1_{ \{\tilde X_{t} \ge nx \}} | \tilde X_{\tau _k} = y ),$$
so
$$ E^0 ( M_{t}  1_{ \{\tilde X_{t} \ge nx \}} ) \leq \sum_{k \, even} g_k \, E^0 ( M_{0,\tau_k}  ) =  \sum_{k \, even}
 g_k .$$

Realize that on $\{\tilde X_{t} \ge nx \} \cap \{ \tau_k \le t <\tau_{k+1}\}$  ($k$ even) we have $\{ \tilde{X}_s \geq nx : s \in [\tau_k, t] \}$.
Recall that $x$ is such that $x_i>0$ for all $i$. Using the assumption on $\phi$ and $\tilde \phi$, we get that asymptotically in $n$:
\begin{align*}
E^0 \Big( 1_{\tau_k \le t <\tau_{k+1} }1_{ \{ \tilde X_{t}  \ge nx \} } & \exp \Big(  \int_{\tau_k}^{t}
\sum_{j=1}^N 1_{ \{ X_{s^-}>0 \}} \log({\tilde \phi_j(X_{s-}) \over \phi_j(X_{s-})})dN^j_s \\
& - \int_{\tau_k}^{t} \sum_{j=1}^N 1_{ \{ X_{s^-} >0 \}} ( \tilde \phi_j(X_{s}) - \phi_j(X_{s}) )ds \Big) \Big) \\
 &\le E^{0}\Big(   1_{\tau_k \le t <\tau_{k+1}}  \exp \{ \sum_{j=1}^N {C_1 \over n} (N^j_t +t) \}  1_{ \{ \tilde X_{t}  \ge nx \} } \Big),\\
 &\le E^{0}\Big(   1_{\tau_k \le t <\tau_{k+1}}  \exp \{ {C_1 N \over n} (\bar N_t +t) \}  1_{ \{ \tilde X_{t}  \ge nx \} } \Big)
\end{align*}
where $\bar{N}_t$ is a Poisson process with parameter $\bar{\lambda}$ equal as the maximum of the parameters of the Poisson processes $\{N^j_t\}_{j=1}^N$.\\

Summing these inequalities and using  H\"{o}lder's inequality 

\begin{align*}
E^0 ( M_{t}  1_{ \{\tilde X_{t} \ge nx \}} ) & \le 
  \exp \{ {C_1 Nt \over n} \} \, \mathbb{P}^{0} ( \tilde X_{t}  \ge nx )^{ \frac{p-1}{p} }  \exp \{ \bar{ \lambda } t [  e^{C_1 pN \over n} - 1] \} ^{\frac{1}{p}}.  
\end{align*}

We now choose $t_n= C_2 n$ such the result of Lemma \ref{lem:contpi} holds, and the sequence $p_n= n$ to obtain
\begin{align*}
{1 \over n }\log \mathbb{P}^0 ( X_{t_n} \geq nx)  
& \le {C_1 \, N C_2  \over n}  + \frac{ \bar{\lambda} C_2  (\exp (C_1 N ) -1 )}{ n }  \\
& +  {n-1 \over  n^2 } \log \mathbb{P}^{0} ( \tilde X_{t_n}  \ge nx ) ,\\
{1 \over n} \log( \pi^{PF} (\{nx\}^{\uparrow})) + { C_3  \over n } \exp ( - C_4 \, n )   & \le  O (n^{-1})  + \Big(\frac{n-1}{n}\Big){1 \over n } \log(  \pi^{mPF} (\{nx\}^{\uparrow}))  \\
& +  { \tilde C_3 (n-1) \over n^2 } \exp ( - \tilde C_4 \, n ).
\end{align*}

So we have
$$ {1 \over n} \log( \pi^{PF} (\{nx\}^{\uparrow})) \leq  {1 \over n} \log( \pi^{mPF}(\{nx\}^{\uparrow})) + O(n^{-1}).$$
We also have the converse of last inequality using the same arguments interchanging $X$ and $\tilde{X}$. \\

Assume now that $x_i>0$ for $i \in \U$, while $x_i=0$ for $i \in \S$, where $\U,\S \subseteq \{1,...,N\}$. 
If classes of  $\U$ are independent of classes $\S$, the result is obvious.
Hence suppose that there exists a coupling between classes of $\U$ and $\S$.
We shall make use of the following fact:
\begin{lemma}
Given the definition of Proportional fairness, if class $i \in \S$ and class $j \in \U$ are coupled, then there exist $K>0$ such that
if $y_i \le n \epsilon$ and $y_j \ge n$ for all $j \in \U$, then:
$$|\phi_i(y) - \tilde \phi_i(y) |\le K \epsilon,$$
\end{lemma}

\bp  

If at least two classes are coupled, there exists a Lagrange multiplier $\alpha>0$ and some positive constants $c_i>0$ and $c_j$ (with at least one $j$ such that $c_j>0$) such that
$${\partial  U^1(y,\eta_i) \over \partial \eta_i} = {x_i \over \eta_i}=\alpha  c_i  , $$
$${\partial  U^1(y,\eta_i) \over \partial \eta_j} = {x_j \over \eta_j}= \alpha c_j . $$
Combined with the fact that there exists $c>0$ such that $0<c \le \sum_{j=1}^N \eta_j \le C$, we obtain that
$$  \tilde c  {y_i \over y_i+ \sum_{j  \in \U_i} y_j} \le \phi_i(y) \le  \tilde C  {y_i \over y_i+ \sum_{j  \in \U_i} y_j} ,$$
where $\U_i \subset \U$. Hence if $y_i \le n \epsilon$ and $y_j \ge n$ for all $j \in \U$
$$\phi_i(y) \le \tilde C { \epsilon \, n \over 1+ n} \le \tilde K \epsilon.$$
Using the control on the modified proportional fair allocation for $x_i \ge 1$:
$$\tilde \phi_i(y) \le \exp\Big({1 \over y_i}\Big) \phi_i(y)\le  e \phi_i(y) \le e \tilde K \epsilon.$$

\ep

Now we do a finer classification of indexes in $\S$. Let $\{ \epsilon_n \}$ be a sequence such that $\epsilon_n$ goes to $0$ as $n$ grows, and
$\{ t_n=C_2 n \}$. For each $n$, we define the sets
$$ \bar \S _n = \{ i \in S :  X_i(t_n) \geq n \epsilon_n \}, \quad \underline{ \S}_n = \{ i \in S :  X_i(t_n) < n \epsilon_n \} .$$
We look at the set of events:
$$\A= \{  X_i(t_n)  \ge n x_i,  i \in \U_R, \, X_i(t_n) \ge n \epsilon _n , i \in \bar{\S}_n, \,  X_i(t_n) \le n \epsilon _n, i \in \underline{\S}_n \}.$$
Recall that
$$ | \log \Big( \frac{ \phi_i(y) }{\tilde \phi_i(y)} \Big) |  \le 1,$$ 
and remark that on $A \in \A$, using the previous lemma, the process counting the number of downwards jumps in direction $j \in \S$ 
which are not common for $X$ and $\tilde X$ is dominated by a Poisson process $N^{\epsilon_n}$ of intensity $C_5 \epsilon_n$  independent of $N^j$, $j \in \S$. \\
Hence:
\begin{eqnarray*}
& & E^0  1_{\tau_k \le t <\tau_{k+1} }1_{ A } M_{\tau_k,t}  \\
\nonumber
 &\leq&  E^{0} \Big( \exp \Big( \sum_{j \in \U} {C_1 \over n} (N^j_t +t) + \sum_{j \in \bar{S}_n} {C_1 \over n \epsilon_n } (N^j_t +t)
 + \sum_{j \in \underline{S}_n }  [ C_5  N^{\epsilon_n} _t +  (C_6 \epsilon_n)t ] \Big)     1_{\tau_k \le t <\tau_{k+1}} 1_{ A } \Big) .
\end{eqnarray*}
and,
\begin{align}\label{ineq:14}
\nonumber
E^0 ( M_{t}  1_{A} )  & \le  E^{0} \Big( \exp \Big( \sum_{j \in \U} {C_1 \over n} (N^j_t +t) + \sum_{j \in \bar{S}_n} {C_1 \over n \epsilon_n } (N^j_t +t)
 + \sum_{j \in \underline{S}_n }  [ C_5  N^{\epsilon_n} _t +  (C_6 \epsilon_n)t ] \Big)  1_{ A } \Big) \\
\nonumber & \le  E^{0} \Big( \exp \Big(  {\# \U C_1 \over n} (\bar{N}_t +t) +   {\# \bar{S}_n C_1 \over n \epsilon_n } (\bar{N}_t +t)
 +  \# \underline{S}_n [  C_5  N^{\epsilon_n} _t +  C_6 \epsilon_n t ] \Big)  1_{ A } \Big) \\
\nonumber &\le  E^{0} \Big(  \exp \{ {C_1 pN \over n} (\bar{N}_t +t) \} \Big)^{\frac{\# \U}{pN}}   E^{0} \Big(  \exp \{ {C_1 pN \over n \epsilon_n} (\bar{N}_t +t) \} \Big)
^{\frac{\# \bar{S}_n }{pN}}  \\
& E^{0} \Big(  \exp \{ pN [C_5  N^{\epsilon_n}_t + C_6 \epsilon_n \, t ] \} \Big)^{\frac{\# \underline{S}_n}{pN}}  \mathbb{P}^{0} ( \tilde X_{t}  \ge nx )^{ \frac{p-1}{p} } ,
\end{align}
where we used H\"older's inequality for $1= \frac{ \# U}{pN} + \frac{\# \bar{S}_n}{pN} + \frac{\# \underline{S}_n}{pN}+ \frac{p-1}{p}$ . \\

We now bound the speed of convergence, in the large deviations scale, of each of the three first factors on the right-hand side of inequality (\ref{ineq:14}). 
Choosing $t_n= C_2 n$, $p_n= C_7 \log n$, $\epsilon_n= n^{-\frac{1}{2}} $, we have:
\begin{align*}
\frac{1}{n} \log E^{0} \Big(  \exp \{ {C_1 N p_n \over n} (\bar{N}_{t_n} +t_n) \} \Big)^{\frac{\# \U}{N p_n }} & = 
\frac{1}{n} \log \Big( \exp \{ {  \# \U C_1 t_n \over  n} \} \exp\{ \bar{\lambda} t_n [e^{\frac{C_1Np_n}{n}}-1]
\}^{\frac{\# \U}{Np_n}} \Big) \\
& = \frac{\# \U C_1 C_2}{n} + \frac{\# \U \bar{\lambda} t_n }{N n p_n} [e^{\frac{C_1Np_n}{n}}-1] \\
&= \frac{\# \U C_1 C_2}{n} + \frac{\# \U \bar{\lambda} C_2 }{N C_7 \log n} [e^{\frac{C_1 C_7 N \log n}{n}}-1] \\
& = O( n^{-1}).
\end{align*}

We do similarly for the second factor:
\begin{align*}
\frac{1}{n} \log E^{0} \Big(  \exp \{ {C_1 N p_n \over \epsilon_n n} (\bar{N}_{t_n} +t_n) \} \Big)^{\frac{\# \bar{ \S}_n}{N p_n }} & = 
\frac{1}{n} \log \Big( \exp \{ {  \# \bar{ \S}_n C_1 t_n \over \epsilon_n n} \} \exp\{ \bar{\lambda} t_n [e^{\frac{C_1Np_n}{ \epsilon_n n}}-1]
\}^{\frac{\# \bar{\S}_n}{Np_n}} \Big) \\
& = \frac{\# \bar{\S}_n C_1 C_2}{\epsilon_n n} + \frac{\# \bar{\S}_n \bar{\lambda} t_n }{N n p_n} [e^{\frac{C_1Np_n}{\epsilon_n n}}-1] \\
& = \frac{\# \bar{\S}_n C_1 C_2}{\epsilon_n n} + \frac{\# \bar{\S}_n \bar{\lambda} C_2 }{N C_7 \log n} [e^{ \frac{C_1 C_7 N \log n}{\epsilon_n n} }-1] \\
& = O( (n \epsilon_n)^{-1}) = O( n^{-\frac{1}{2}}).
\end{align*}

Finally, for the third factor we have:
\begin{align*}
\frac{1}{n} \log E^{0} \Big(  \exp \{ N p_n [C_5  N^{\epsilon_n}_{t_n} + C_6 \epsilon_n \, t_n ] \} \Big)^{\frac{\# \underline{S}_n}{N p_n}}  &=
\frac{1}{n} \log \Big( \exp \{ C_6 \# \underline{S}_n \epsilon_n t_n  \} \exp \{ \epsilon_n t_n [  e^{C_5N p_n} - 1] \}^{\frac{\# \underline{S}_n}{Np_n}} \Big) \\
&= \frac{  C_6 \epsilon _n \, \# \underline{S}_n t_n}{n} + \frac{\# \underline{S}_n \epsilon_n \, t_n}{N p_n  n}  [  e^{C_5N p_n} - 1]  \\
 & = C_6 C_2 \# \underline{S}_n \epsilon _n \,  + \frac{C_2 \# \underline{S}_n \epsilon_n }{C_7 N \log n }  [  e^{C_5 C_7 N \log n} - 1]  \\
 & = C_6 C_2 \# \underline{S}_n \epsilon _n \,  + \frac{C_2 \# \underline{S}_n \epsilon_n [ n^{C_5 C_7 N } - 1] }{C_7 N \log n }  \\
 & = O(\epsilon_n) + O( \frac{ \epsilon_n n^\epsilon}{\log n} ) = O(n^{- \frac{1}{2} + \epsilon}).
\end{align*}
where we chose $C_7>0$ small enough such that $C_7 C_5 N < \epsilon$, for any given $\epsilon>0$. \\
 
Then along the same lines as in the case $x$ has all entries positive, it follows
$$ {1 \over n} \log( \pi^{PF} (\{nx\}^{\uparrow})) \leq  {1 \over n} \log( \pi^{mPF}(\{nx\}^{\uparrow})) + O( n^{-\frac{1}{2}+\epsilon} ) \quad  \forall \epsilon >0 .$$

\ep

\section{Large deviations for monotone networks with general service time distributions}
\label{sec:MN}

In this Section, we consider a processor sharing network (i.e. a set of processor sharing nodes)
 with a proportionally fair bandwidth allocation. This means that a flow is served at node $i$
with speed ${\phi^{PF}_i(X_t) \over X_i(t)}$. We need the notion of strong monotonicity introduced for instance in \cite{stobounds}.

\begin{definition}
$\phi$ is strongly decreasing if the function $\psi_i$ defined by $\psi_i(x)= { \phi_i(x) \over x_i}$, if $x_i >0$ and $0$ otherwise, is decreasing in $x_j$ for all $j \neq i$.
\end{definition}

 We assume in this Section that the network is monotone, i.e. that $\phi$ is strongly increasing.
 This is verified on  all tree topologies for instance and for many wireless networks instances.
In the proof of the following theorem, we use the monotony of the network to obtain stochastic comparisons. We recall the following Proposition:

\begin{proposition}
Let $X_t$ and $\tilde X_t$ two processes associated with the allocations $\phi$ and $\tilde \phi$.
Suppose that $\phi$ is strongly decreasing, and that ${\phi_i(x) \over x_i}\ge { \tilde \phi_i(x) \over x_i}$ for each $x$. Then for all $t$ and for all service time distribution:
$$ X_t \prec_{st} \tilde X_t .$$
\end{proposition}

\bp

Remark that the processes $X_t$ and $\tilde X_t$ are not Markov in general.
We can however easily construct a sample-path comparison between the processes.
Starting with ordered initial configuration of flows $x \prec y$, the rates at which each flow is served in the network verify:
$${\phi_i(x) \over x_i}\ge{\phi_i(y) \over y_i} \ge { \tilde \phi_i(y) \over y_i}.$$
Hence, a coupling can be constructed such that  $X_t \prec \tilde X_t$ almost surely.

\ep

\

Recall that we define $P^{PF}$ as $P^{PF}(x)=\max_{\eta \in \log (\C) } \langle \eta,x \rangle$,
and that $\log(\phi_i(x))={\nabla}P^{PF}(x).$

\

Now define the discrete potential functions $\overline \Psi$ and $\underline \Psi$ by the recursive formula:

$$\overline \Psi(x)=\max_{i=1\ldots N}(\overline \Psi(x-e_i) -\log(\phi_{i}(x))),$$
$${\underline \Psi}(x)=\min_{i=1\ldots N} \{ \underline \Psi(x-e_i) -\log(\phi_{i}(x)) \}.$$
We call supPF and infPF the reversible allocations associated with the discrete potentials $\overline \Psi$ and $\underline \Psi$ and
define the rate function $R(x)= \langle \log(\lambda),x \rangle -P(x)$. We can now state the main result of this Section:

\begin{theorem}
Suppose the network monotone, i.e. $\phi$ is strongly monotone.
If $\lambda$ is in the interior of  the capacity set $\C$, then
 the allocations Proportional fairness, supPF and infPF as well as Balanced fairness are stable and admit the same large deviation characteristics with rate function $R$.

\end{theorem}

\pagebreak

\bp

Given the definitions of the discrete potentials $\overline \Psi$ and $\underline \Psi$, observe that
 there exists two paths $\overline \P(x)$ and $\underline \P(x)$ from $0$ to $x$ such that:
$$\overline \Psi(x)=-\sum_{k : (i_k,z_k) \in \overline \P(x)} \log(\phi_{i_k}(z_k)),$$
$$\underline \Psi(x)=-\sum_{k : (i_k,z_k)  \in \underline \P(x)} \log(\phi_{i_k}(z_k)),$$
where with a slight abuse of notations the indexes $i(z)$ correspond to the indexes
defined by the specific paths $ \overline \P(x)$ and  $\underline \P(x)$. \\

Let $\overline X$ and $\underline X$ the process associated with the balance function $\overline \Psi$ and $\underline \Psi$.
Recall now that $P$ is the potential associated with the proportional fair allocation, and since we can define this potential up to an additive constant, let us choose it such that $P(0)=0$.
So, we can write that for any path $\P(x)$ going from $0$ to $x$, $P(x)=P(x)-P(0)=\sum_{k : (i_k,z_k)  \in \P(x)} \mathbf{D}_{i_k} P(z_k)$.
Choosing $\P(x)=\bar \P(x)$:
\begin{eqnarray}\label{eq:1}
\nonumber
|\overline \Psi(n x)-P(n x)|&=& | \sum_{ k : (i_k,z_k)  \in \bar \P(n x)} \mathbf{D}_{i_k} P(z_k)- \log( \phi_{i_k}(z_k) |,\\
&\le&  \sum_{ k : (i_k,z_k)  \in \bar \P(n x)} | \mathbf{D}_{i_k} P(z_k)- \log( \phi_{i_k}(z_k) |.
\end{eqnarray}
 Recalling that for $x_i\ge 1$:
 $$|P(x)-P(x-e_i)-\log(\phi_{i}(x))| \le {1 \over x_i },$$ 
 we can bound (\ref{eq:1}) by:
\begin{eqnarray}\label{eq:2}
\nonumber |\overline \Psi(n x)-P(n x)| &\le& \sum_{(i_k, z_k) \in \bar \P(n x)} { 1 \over z_{i_k}},\\
&\le&  C  \sum_{i=1}^{n} {1 \over i} \le  C' \log(n).
\end{eqnarray}

Define the invariant measures (not necessarily stationary) $\overline \pi$ and $\underline \pi$ of the processes $\overline X \prec X \prec \underline X$.
Since we have constructed these processes from a discrete potential, we have that (see Proposition \ref{prop:rev}):
$$\overline \pi(x)=\lambda^x \exp(-\overline \Psi(x)),$$
$$\underline \pi(x)=\lambda^x \exp(-\underline \Psi(x)).$$
Now observe that:
$${1 \over n}\log(\overline \pi(nx))=  \langle \log(\lambda),x \rangle - {1 \over n}\overline \Psi(nx),$$
$${1 \over n}\log(\underline \pi(nx))=  \langle \log(\lambda),x \rangle - {1 \over n}\underline \Psi(nx).$$
Using (\ref{eq:2}) we obtain that 
$${1 \over n}\log(\overline \pi(nx))= {1 \over n} \langle \log(\lambda),nx \rangle - {1 \over n} P(nx) + O(n^{-1}),$$
$${1 \over n}\log(\underline \pi(nx))= {1 \over n} \langle \log(\lambda),nx \rangle - {1 \over n} P(nx) + O(n^{-1}),$$
and hence, by the $1$-homogeneity of $R$,
$$\lim_{n \to \infty}{1 \over n}\log(\overline \pi(nx))= \lim_{n \to \infty}{1 \over n}\log(\underline \pi(nx))=R(x).$$

Assume $\lambda$ is in the interior of the capacity set $\C$. Using the definition of $P$, this implies that there exists $a>0$
such that
$$R(x) = \langle \log(\lambda),x \rangle -P(x)= \, \langle \log(\lambda)-\log(\phi(x)),x \rangle  \le - a |x|,$$
which implies that the invariant measure $\overline \pi$ and $\underline \pi$ are summable.
This hence proves that the stationary distributions of both processes $\overline X$ and $\underline X$ are well defined for Markovian dynamics while the reversibility condition (i.e. the balance property of the service rates) implies the
insensitivity of the stationary distribution to the service time distribution \cite{perf,zachary2007}. Hence, for any service time distribution, the networks with allocation infPF and supPF
admit the stationary distribution $ C_1\overline \pi$ and $C_2 \underline \pi$ where $C_1$ and $C_2$ are normalizing constants.
Now using the assumption of monotonicity of the network, we obtain that:
\begin{equation}
P(|\overline X| \ge n)\ge P(| X| \ge n) \ge P(|\underline X| \ge n),
\end{equation}
and by Cram\'er's theorem we conclude:
$$ \lim_{n \rightarrow \infty }{1 \over n} \log P(|X|\ge n)= \lim_{n \rightarrow \infty}{1 \over n} \log P(|X|= n).$$
It has further been proven in \cite{massoulie2007} that mPF and BF have $R(x)$ as rate functions.

\ep

\section{Conclusion}
We proved that the stationary measure of the number of flows in progress in a bandwidth sharing network functioning under the proportional fair allocation
shares the same large deviations characteristics with the stationary meausure of the number of flows in progress of the same network under the balanced fair allocation. 
This formalizes the idea that in long excursions proportional fair allocation behaves similarly to the most efficient insensitive allocation.

\

\bibliographystyle{abbrv}
\bibliography{LD2010}

\end{document}